**Freight Traffic Assignment Methodology for Large-Scale Road-Rail Intermodal Networks**

Md Majbah Uddin
Department of Civil and Environmental Engineering
University of South Carolina
300 Main Street, Columbia, SC 29208
Email: muddin@email.sc.edu

and

Nathan Huynh (*Corresponding Author*)
Department of Civil and Environmental Engineering
University of South Carolina
300 Main Street, Columbia, SC 29208
Email: huynhn@cec.sc.edu

Length of paper: Text = 6,320 words; Number of Figures and Tables = 6; Total = 7,820 words

**ABSTRACT**

This paper proposes a methodology for freight traffic assignment in large-scale road-rail intermodal networks. To obtain the user-equilibrium freight flows, a path-based assignment algorithm (gradient projection) is proposed. The developed methodology is tested on the U.S. intermodal network using the 2007 freight demands for truck, rail, and road-rail intermodal from the Freight Analysis Framework, version 3, (FAF3). The results indicate that the proposed methodology's projected flow pattern is similar to the FAF3 assignment. The proposed methodology could be used by transportation planners and decision makers to forecast freight flows and to evaluate strategic network expansion options.

**INTRODUCTION**

Freight transportation is a vital component of the U.S. economy. Its chief role is to move raw materials and products in an efficient manner (*1*). The U.S. has the largest freight transportation system in the world (*2*). It moved, on average, 54 million tons worth nearly $48 billion of freight each day in 2012 and the majority of freight was transported by either truck or rail (67% by truck and 10% by rail). The freight volume is expected to increase to 78 million tons (about 45%) by the year 2040 (*3*). In recent years, intermodal transport is becoming an increasingly attractive alternative to shippers, and this trend is likely to continue as governmental agencies are considering policies to induce a freight modal shift from road to intermodal to alleviate highway congestion and emissions. Moreover, significant social benefits such as enhanced highway safety, reduction in need for building highways, etc. can be obtained by greater use of intermodal (*4*). For an intermodal shipment, the long-haul is carried out by rail, inland waterway or sea, and the initial and/or final short-haul is carried out by road (*5*). With the growth of intermodal transportation, there is a need by transportation planners and decision makers to forecast freight flows on the intermodal networks and to evaluate strategic network expansion options. Well-informed infrastructure, economic, and environmental planning depends on effective freight forecasting (*6*) which is obtained from the freight assignment step. The multimodal nature of the freight movement presents an additional layer of complexity to the freight assignment problem. Additionally, freight demand and cost data are not as readily available. To this end, this paper proposes an integrated freight assignment methodology that considers road, rail and intermodal shipments.

The assignment of freight over multimodal networks has been studied by many researchers in the past few decades. Crainic et al. (*7*) developed a nonlinear optimization model to route freight train, schedule train services and allocate classification work between yards. Guelat et al. (*8*) proposed a Gauss-Seidal-Linear approximation algorithm to assign multiproduct in a multimode network for strategic planning. Their algorithm was implemented in a strategic analysis tool named "strategic transportation analysis (STAN)" and solved a system-optimal problem with the objective of minimizing the total cost at arcs and node transfer. Their solution algorithm considered intermodal transfer costs in the computation of shortest paths. Chow et al. (*6*) considered nonlinear inverse optimization for the freight assignment using a variant of STAN, and calibrated their model to work for both user-equilibrium and system-optimal conditions. The freight network equilibrium model (FNEM) developed by Friesz et al. (*9*) considered the combined role of shipper-carrier. Using the shipper and carrier sub-models FNEM provided the route choice decisions for both shippers and carriers on a multimodal freight network with nonlinear cost and delay function. By solving a variational inequality (VI) problem on the railway network Fernández et al. (*10*) developed a strategic railway freight assignment model. Agrawal and Ziliaskopoulos (*11*) also used the VI approach for freight assignment to achieve market equilibrium where no shipper can reduce its cost by changing carrier. In their model, shippers were assumed to have user-equilibrium behavior with the objective of minimizing cost without any consideration about other shippers in the market, whereas carriers followed a system-optimal behavior with the objective of optimizing their system (complete operation).

Loureiro and Ralston (*12*) proposed a multi-commodity multimodal network design model to use as a strategic planning tool; the model assumed that the goods are shipped at minimum total generalized cost and used path-based user-equilibrium assignment algorithm to assign freight flows over the network. Kornhauser and Bodden (*13*) analyzed highway and intermodal railway-highway freight network by routing freight over the network using a minimum cost path-finding algorithm and presented results as density map. Arnold et al. (*14*) proposed a modeling framework for road-rail intermodal network, but the main purpose of their model was to optimally locate intermodal terminals by minimizing transportation cost of shipments. Mahmassani et al. (*15*) developed a dynamic freight network simulation-assignment model for the analysis of multiproduct intermodal freight transportation systems. The intermodal shortest path was calculated based on the link travel costs and node transfer delays. Zhang et al. (*16*) validated the Mahmassani et al. model by applying it to a Pan-European rail network. Using a bi-level programming, where lower-level problem finds the multimodal multiclass user traffic assignment and upper-level problem determines the maximum benefit-cost ratio yielding network improvement actions, Yamada et al.

(*17*) developed a multimodal freight network model for strategic transportation planning. Chang (*18*) formulated a route selection problem for international intermodal shipments considering multimodal multi-commodity flow. The model was formulated to consider multiple objectives, scheduled modes and demanded delivery times, and economies of scale. Hwang and Ouyang (*19*) used the user-equilibrium approach to assign freight shipments onto rail networks which were represented as directed graphs.

Based on the above review, to date, no model has been developed to comprehensively assign freight flows that are transported via multiple modes (road-only, rail-only, and road-rail intermodal) under equilibrium conditions. This paper seeks to fill this gap in the literature by developing such a model. Specifically, given a set of freight demands between origins and destinations and designated modes (road-only, rail-only, and intermodal), the model seeks an equilibrium assignment that minimizes the total transportation cost (i.e., travel time) for the freight transport network. To solve the proposed model, a path-based algorithm, based on the gradient projection algorithm proposed by Jayakrishnan et al. (*20*), is adopted. The gradient projection algorithm is chosen because it has been shown to converge faster than the conventional Frank-Wolfe algorithm (*21*) and outperform other path-based algorithms (*22*).

To model congestion effects in a network at the planning level, link performance functions are often used, which express the travel time on a link as a function of link flow. For highways, the standard Bureau of Public Road (BPR) link performance function is commonly used. For rail, a few functions have been proposed (*19, 23, 24*). Recently, Borndörfer et al. (*25*) suggested a link performance function for freight rail network; its functional form is similar to its highway counterpart. When applying these types of functions, it is necessary to calibrate the parameters to capture local and regional effects. In this study, the function proposed by Borndörfer et al. is adopted and calibrated to reflect characteristics of the U.S. rail infrastructure.

To validate the proposed model, the projected equilibrium freight flow pattern on the U.S. intermodal network is compared against the Freight Analysis Framework, version 3, (FAF3) network flow assignment pattern. FAF3 is the most comprehensive public source of freight data in the U.S. The FAF3 database contains aggregated freight data for 131 origins, 131 destinations, 43 commodity classes, and 8 modal categories (*26*). It should be noted that the FAF3 flow values are not absolute. Rather, the FAF3 flows are estimated using models that disaggregate interregional flows into flows between localities (based on geographic distributions of economic activity) and then these flows are assigned to individual highways using average payloads per truck to produce truck counts. Thus, the FAF3 flow values could be different from actual truck counts.

The remainder of the paper is organized as follows. The formulation of the model and solution algorithm is presented in the next section. Then, the application of the model is demonstrated using the U.S. intermodal network. The study's conclusions are given in the last section.

## MODELING AND ALGORITHMIC FRAMEWORK

This study takes a system's view and assumes that in the long run the activities carried out by shippers and carriers will lead to equilibrium where the cost of any shipment cannot be lowered by changing mode and/or route. The freight logistics problem has two levels. The first and upper level involves decisions by shippers in selecting a carrier, and the second and lower level involves decisions by the carriers in minimizing the shipment times. The modeling framework proposed here (i.e. freight traffic assignment) is for the lower level. Therefore, it is assumed, the cost on all used paths via different modes (road-only, rail-only, and intermodal) is equal for each O-D demand pair and equal to or less than the cost on any unused path at equilibrium (*27*).

### Notation

| | |
|---|---|
| $N$ | Set of nodes in the network |
| $A$ | Set of links in the network |
| $N_c$ | Set of freight zone centroid nodes in the network |
| $N_t$ | Set of road nodes in the network |

$N_l$       Set of rail nodes in the network

$A_t$       Set of road links in the network

$A_l$       Set of rail links in the network

$A_f$       Set of terminal links in the network

$R$       Set of origins in the network, $R \subseteq N$

$S$       Set of destinations in the network, $S \subseteq N$

$r$       Origin zone index, $r \in R$

$s$       Destination zone index, $s \in S$

$x_a$       Flow on link $a$, $a \in A$

$t_a(\omega)$       Travel time on link $a$ for a flow of $\omega$

$f_k^{rs}$       Flow on path $k$ connecting $r$ and $s$

$f_{\overline{k}^n}^{rs}$       Flow on shortest path connecting $r$ and $s$

$q_t^{rs}$       Freight truck demand from $r$ to $s$

$q_l^{rs}$       Freight train demand from $r$ to $s$

$q_i^{rs}$       Freight intermodal demand from $r$ to $s$

$K_t^{rs}$       Set of paths with positive truck flow from $r$ to $s$

$K_l^{rs}$       Set of paths with positive train flow from $r$ to $s$

$K_i^{rs}$       Set of paths with positive intermodal flow from $r$ to $s$

$T$       Set of available terminals for transfer of shipments

## Formulation

Consider a network which is represented by a directed graph $G = (N, A)$, where $N$ is the set of nodal points of the network ($N = N_c \cup N_t \cup N_l$), while $A$ is the set of links joining them in the network ($A = A_t \cup A_l \cup A_f$). In the network, nodal points are made of three node sets: zone centroid represented by nodes ($N_c$), road intersections ($N_t$), and rail junctions ($N_l$). On the other hand, network links are formed by three sets: road segments ($A_t$), rail tracks ($A_l$), and terminal transfer links ($A_f$). Note that road-rail intermodal terminals are modeled as links and that flows are bi-directional on these links. Furthermore, their end nodes have different modes (one from the set $N_t$ and the other from the set $N_l$). For truck traffic demand $q_t^{rs}$ from origin $r \in R$ to destination $s \in S$ and a set of paths that connect $r$ to $s$ for each O-D pair $K_t^{rs}$, the independent variables are a set of path flows $f_k^{rs}$ that satisfy the demand ($\sum_{k \in K_t^{rs}} f_k^{rs} = q_t^{rs}$). Similarly, the path flows for train and intermodal on path-sets, $K_l^{rs}$ and $K_i^{rs}$, satisfy their respective demands ($q_l^{rs}$ and $q_i^{rs}$) from $r$ to $s$. Note that the path-set for intermodal consists of paths formed by links from both road and rail segments of the network. Therefore, the total freight flow on a road segment ($a \in A_t$) is the sum of the road-only flows and intermodal flows. Similarly, the total freight flow on a rail segment ($a \in A_l$) is the sum of the rail-only and intermodal flows. The user-equilibrium model for this problem is formulated as follows.

$$\text{Min } Z = \sum_{a \in A_t} \int_0^{x_a} t_a(\omega)\,d\omega + \sum_{a \in A_l} \int_0^{x_a} t_a(\omega)\,d\omega \tag{1}$$

$$\text{Subject to } \sum_{k \in K_t^{rs}} f_k^{rs} = q_t^{rs}, \qquad \forall r \in R, s \in S \tag{2}$$

$$\sum_{k \in K_r^{rs}} f_k^{rs} = q_l^{rs}, \qquad \forall r \in R, s \in S \tag{3}$$

$$\sum_{k \in K_i^{rs}} f_k^{rs} = q_i^{rs}, \qquad \forall r \in R, s \in S \tag{4}$$

$$x_a = \sum_{r \in R} \sum_{s \in S} \sum_{k \in K_t^{rs}} f_k^{rs} \delta_{ka}^{rs} + \sum_{r \in R} \sum_{s \in S} \sum_{k \in K_i^{rs}} f_k^{rs} \delta_{ka}^{rs}, \quad \forall a \in A_t \tag{5}$$

$$x_a = \sum_{r \in R} \sum_{s \in S} \sum_{k \in K_l^{rs}} f_k^{rs} \delta_{ka}^{rs} + \sum_{r \in R} \sum_{s \in S} \sum_{k \in K_i^{rs}} f_k^{rs} \delta_{ka}^{rs}, \quad \forall a \in A_l \tag{6}$$

$$f_k^{rs} \geq 0, \qquad \forall k \in K_t^{rs}, k \in K_l^{rs}, k \in K_i^{rs}, r \in R, s \in S \tag{7}$$

where,

$$\delta_{ka}^{rs} = \begin{cases} 1 & \text{if link } a \text{ is on path } k \text{ connecting } r \text{ and } s \\ 0 & \text{otherwise} \end{cases}$$

The objective function (1) states that the total travel time for both segments (road and rail) associated with the flows between origins and destinations are to be minimized. Constraints (2) - (4) ensure that all freight demands are assigned to the network. Constraints (5) and (6) are definitional constraints that compute link flows. Lastly, constraint (7) ensures non-negative flows.

To model congestion effects in a network at the planning level, link performance functions are often used, which express the travel time on a link as a function of link flow. For highways, the standard Bureau of Public Road (BPR) link performance function, named after the agency which developed it, is commonly used. For rail, a few functions have been proposed (*19, 23, 24*). Recently, Borndörfer et al. (*25*) suggested a link performance function for freight rail network; its functional form is similar to its highway counterpart. When applying these types of functions, it is necessary to calibrate the parameters to capture local and regional effects. In this study, the function proposed by Borndörfer et al. is adopted and calibrated to reflect characteristics of the U.S. rail infrastructure. The link performance functions have the following form.

$$t_a(x_a) = t_{o,t}\left(1 + 0.15\left(\frac{x_a}{C_a}\right)^4\right), \qquad \forall a \in A_t \tag{8}$$

$$t_a(x_a) = t_{o,l}\left(1 + \left(\frac{x_a}{C_a}\right)^\beta\right), \qquad \forall a \in A_l \tag{9}$$

where $t_{o,t}$ and $t_{o,l}$ are the free-flow travel time for road and rail links, respectively and $C_a$ is the capacity of the link. In equation (9), $\beta$ represents the penalty rate and its value can be 2, 4, 7, 15 (*25*). In this paper $\beta$ is calibrated to capture characteristics of the rail segment of the U.S. intermodal network. Calibration involved changing the value of $\beta$ such that the computed train delay resulted in realistic flow pattern. The functional form of equation (9) indicates that the travel time on rail links is more sensitive to flow when it is near capacity than that of road links.

Figure 1 illustrates the methodology used to calculate the intermodal shortest path. Part (a) of Figure 1 shows the typical intermodal freight transport elements that are used to ship goods from an origin to a destination; a typical shipment would go through two intermodal terminals. Part (b) of Figure 1 shows the corresponding network structure. The intermodal path is made up of the node sequence: $b \rightarrow c \rightarrow d \rightarrow e \rightarrow f \rightarrow g$. Thus, given $b$ and $g$, the objective of the shortest path algorithm is to find nodes $c$, $d$, $e$, and $f$ that result in the least travel time. Delays are incurred at intermodal terminals due to the transfer of modes and storage (when the timing of inbound and outbound trains does not coincide). This terminal delay is considered as terminal link delay $(t_a, \forall a \in A_f)$ in the path travel time calculation.

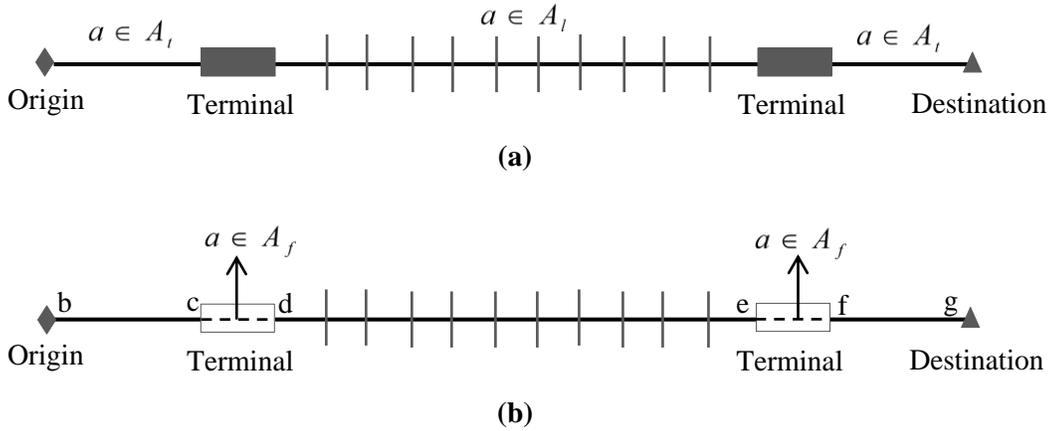

**(a)**

**(b)**

**FIGURE 1 Shortest path calculation considering terminal: (a) basic intermodal structure and (b) modeled structure.**

## Solution Algorithm

A path-based algorithm (gradient projection) is used to solve the proposed user-equilibrium assignment problem. The adopted gradient projection algorithm is based on the Goldstein-Levitin-Polyak gradient projection method formulated by Bertsekas (*28*) and modified by Jayakrishnan et al. (*20*) to solve the traffic assignment problem. In this study, this algorithm is further modified to address the assignment of freight demands that can be transported via three different modes: road-only, rail-only, and intermodal. Additionally, the algorithm is modified to consider intermodal terminals in the network. The iterative steps of the algorithm are as follows.

**Step 0:** Initialization.

Set $t_a = t_a(0), \forall a \in A$ and select terminals for all O-D pairs. Assign O-D demands $q_t^{rs}$, $q_l^{rs}$, and $q_i^{rs}$ on the shortest path calculated based on $t_a, \forall a \in A_t$, $t_a, \forall a \in A_l$, and $t_a, \forall a \in A$, respectively and initialize the path-sets $K_t^{rs}$, $K_l^{rs}$, and $K_i^{rs}$ with the corresponding shortest path for each O-D pair $(r, s)$. This yields path flows and link flows. Set iteration counter $n = 1$.

**Step 1:** For each O-D pair $(r, s)$:

**Step 1.1:** Update.

Set $t_a(n) = t_a(x_a(n)), \forall a \in A$. Update the first derivative lengths (i.e., path travel times at current flow): $d_{kt}^{rs}(n), \forall k \in K_t^{rs}$, $d_{kl}^{rs}(n), \forall k \in K_l^{rs}$, and $d_{ki}^{rs}(n), \forall k \in K_i^{rs}$.

**Step 1.2:** Direction finding.

Find the shortest path $\bar{k}_t^{rs}(n)$ based on $t_a(n), \forall a \in A_t$. If different from all the paths in $K_t^{rs}$, add it to $K_t^{rs}$ and record $d_{\bar{k}_t^{rs}(n)}^{rs}$. If not, tag the shortest among the paths in $K_t^{rs}$ as $\bar{k}_t^{rs}(n)$.

Repeat this procedure for $K_l^{rs}$ and $K_i^{rs}$ to find $d_{\bar{k}_l^{rs}(n)}^{rs}$ and $d_{\bar{k}_i^{rs}(n)}^{rs}$ based on $t_a(n)$, $\forall a \in A_l$ and

$t_a(n)$, $\forall a \in A$, respectively.

**Step 1.3:** Move.

Set the new path flows for $K_t^{rs}$.

$$f_k^{rs}(n+1) = \max\left\{0, f_k^{rs}(n) - \frac{\alpha(n)}{s_k^{rs}(n)}\left(d_{kt}^{rs}(n) - d_{\bar{k}_t^{rs}(n)}^{rs}\right)\right\}, \quad \forall\; k \in K_t^{rs}, k \neq \bar{k}_t^{rs}$$

where,

$$s_k^{rs}(n) = \sum_a \frac{\partial t_a^{rs}(n)}{\partial x_a^{rs}(n)}, \qquad\qquad\qquad \forall\; k \in K_t^{rs}$$

$a$ denotes links that are on either $k$ or $\bar{k}_t^{rs}$, but not on both. $\alpha(n)$ is the step-size (value = 1 (*20*))

Also, $f_{\bar{k}^n_t}^{rs}(n+1) = q_t^{rs} - \sum f_k^{rs}(n+1)$, $\qquad \forall k \in K_t^{rs}, k \neq \bar{k}_t^{rs}(n)$

Follow this procedure to find new path flows for $K_l^{rs}$ and $K_i^{rs}$.

From path flows find the link flows $x_a(n+1)$.

**Step 2:** Convergence test.

If the convergence criterion is met, stop. Else, set $n = n+1$ and go to step 1.

For rail networks, the same infrastructure (i.e., rail tracks) is often shared by traffic flow in both directions. To model this feature, two separate directed links in opposite directions are used instead of one bi-directional link. These two links share the same properties such as length and capacity. Moreover, the link delay on any one link is dependent on the flow on it, as well as the flow on the opposite link (see (*19*) for details). Due to the use of this modeling method, the link performance function shown in equation (9) needs to be modified. The modified version is shown in equation (10), where $x_a$ is the link flow from node $i$ to node $j$ and $x_{a'}$ is the flow from node $j$ to node $i$. Equation (1) also needs to be modified and its modified version is shown in equation (11). The rest of the model is the same and the above solution algorithm remains applicable for solving the modified model.

$$t_a(x_a + x_{a'}) = t_{o,l}\left(1 + \left(\frac{x_a + x_{a'}}{C_a}\right)^\beta\right), \qquad \forall a \in A_l \tag{10}$$

$$Z = \sum_{a \in A_r} \int_0^{x_a} t_a(\omega)\,d\omega + \sum_{a \in A_l} \int_0^{x_a + x_{a'}} t_a(\omega)\,d\omega \tag{11}$$

The proposed model provides a general framework for addressing different types of freight transport networks and situations. While the highway mode generally allows truck to provide door-to-door service, there may be some situations where trucks are not allowed to traverse certain segments in the network. Similarly, certain rail track segments may be accessible or available to shippers. The proposed model can address this by restricting those links in shortest path calculation, and thus, those restricted links are not considered in the assignment process. The model can also address the situations when some intermodal terminals are not available for routing shipments between certain O-D demand pairs. This can be done by excluding those terminals from the set (*T*) for an O-D demand pair during terminal selection (i.e., initialization step of solution algorithm).

**Special Case (Intermodal Demand Only)**
The proposed model is also applicable for intermodal freight demand assignment, with a few modifications. Given all the network elements and demand ($q_i^{rs}$), the intermodal assignment problem is as follows.

$$\text{Min } Z = \sum_{a \in A_r} \int_0^{x_a} t_a(\omega)\,d\omega + \sum_{a \in A_l} \int_0^{x_a} t_a(\omega)\,d\omega \tag{12}$$

$$\text{Subject to } \sum_{k \in K_i^{rs}} f_k^{rs} = q_i^{rs}, \qquad \forall\, r \in R, s \in S \tag{13}$$

$$x_a = \sum_{r \in R} \sum_{s \in S} \sum_{k \in K_i^{rs}} f_k^{rs} \delta_{ka}^{rs}, \quad \forall\, a \in A \tag{14}$$

$$f_k^{rs} \geq 0, \qquad \forall\, k \in K_i^{rs}, r \in R, s \in S \tag{15}$$

The solution algorithm described previously is also applicable for solving problem (12) - (15). However, path-set $K_i^{rs}$ and shipment demand $q_i^{rs}$ should be considered in the solution algorithm instead of three path-sets and three demands.

**APPLICATION**
To demonstrate the validity of the proposed methodology, the model is applied to the U.S. intermodal network created by Oak Ridge National Laboratory (*29*). Without loss of generality, the network is modified to retain only the primary elements of the network. The assignment problem is investigated from a strategic perspective. Thus, freight flows are assigned to the entire freight transport network without considering any restrictions on highway links, rail links, and intermodal terminals.

**Network Description**
The intermodal network considered is shown in Figure 2. Part (a) shows the detailed version, and part (b) shows the simplified version. As shown in the figure, the intermodal network comprises the U.S. interstates, Class I railroads and road-rail terminals. The squares represent freight zone centroids. The circles represent road-rail terminals. The black lines represent interstates, and the gray lines represent Class I railroads. The simplified network has a total of 1532 links and 301 nodes. The nodes include 120 centroids, 97 road intersections, and 84 rail junctions.

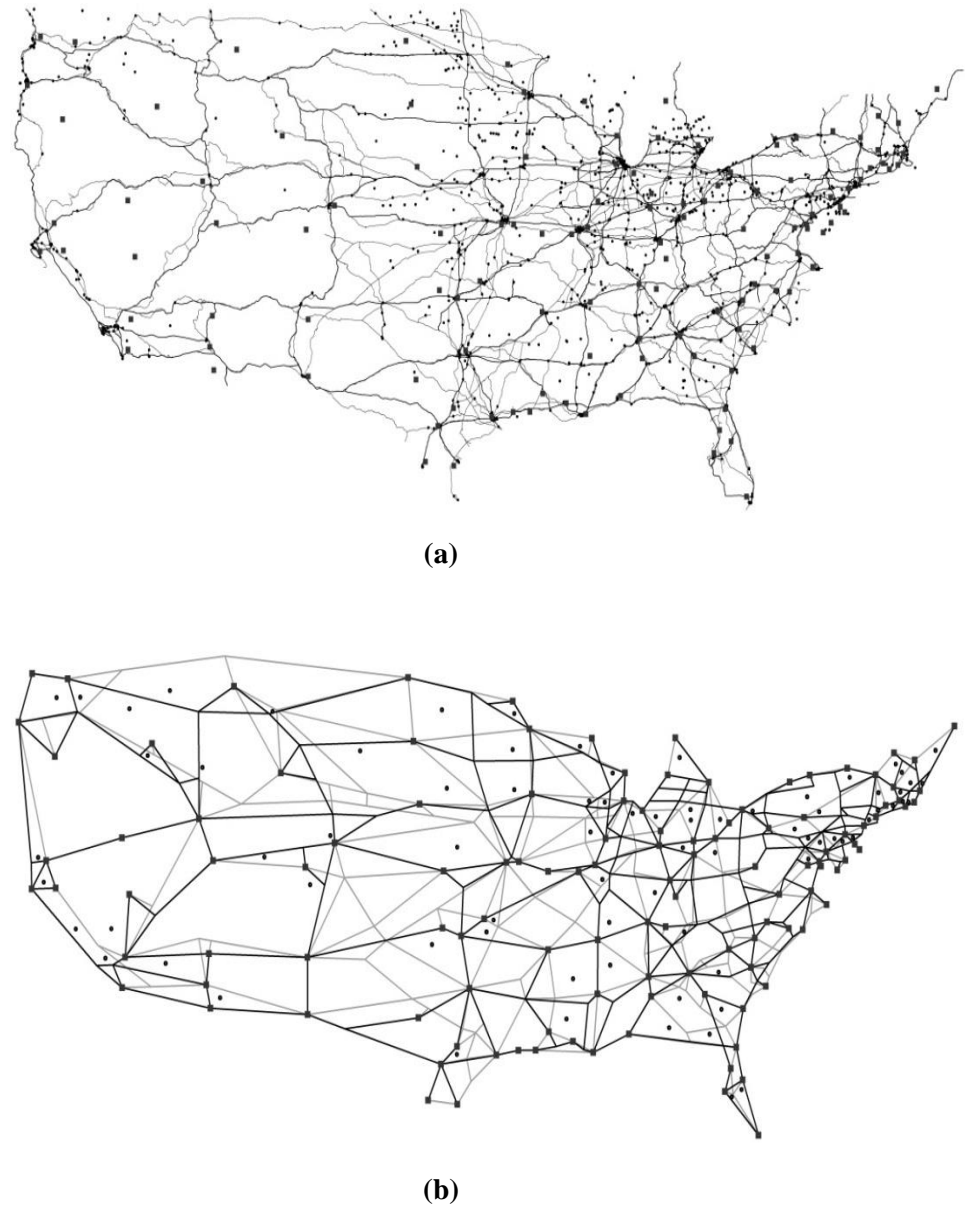

**(a)**

**(b)**

**FIGURE 2 Road-rail transportation networks in the contiguous U.S.: (a) detailed network and (b) simplified network.**

Attributes of the network elements include link lengths, number of tracks, type of control for rail links, etc. The free-flow speed for the road links is calculated using the equation provided in the NCHRP Report 387 (*30*) which requires speed limit as an input. For the rail links, the maximum speed for freight train is taken as 60 mph (*23*). Free-flow travel times for links are calculated using free-flow speeds. Capacities for the rail links are obtained using the number of tracks and type of control for corresponding rail links (*31*). For rural interstates and urban interstates, a capacity of 21,000 veh/lane/day and 19,500 veh/lane/day is used, respectively (*32*). Rail links are assumed to have full capacity, whereas road links are assumed to have reduced capacity due to congestion. In the network considered, contiguous U.S., the total number of freight zones is 120, and hence it is assumed that there are 14,400 possible O-D demand pairs in the network. The freight demands for all O-D pairs are obtained from the FAF3 database (*33*).

The FAF3 procedure to convert tonnage to truck counts (*34*) is used in this study and the key steps are summarized here: (i) compute distance between origin and destination centroid, (ii) using truck allocation factors based on five distance ranges allocate tonnage to five truck types, (iii) convert tonnage assigned to each truck type into their equivalent annual truck traffic values using the truck equivalency factors, which is based on 9 truck body types, (iv) find empty trips using empty truck factors and add empty trips to the loaded trips, (v) aggregate the total annual truck traffic for all body styles together for each truck types, and (vi) sum the traffic for all the truck types. The output of this conversion process is the overall annual truck traffic between the origin and destination. This procedure is carried out for all the demands that are transported by trucks.

The procedure to convert tonnage to trainloads developed by Hwang (*35*) is used in this study. The conversion steps are: (i) group FAF commodity types into 10 types based on similarities, (ii) convert tonnage into equivalent trainloads using average loading weight factors for each commodity group, and (iii) sum the trainloads for all commodity groups. This procedure is carried out for all the demands that are transported by rail.

FAF3 does not provide intermodal demand directly. Thus, to obtain this information the demand recorded as being transported by "multiple modes and mail" is used. To estimate the intermodal demand from this source, several filters are applied. The data are filtered to include only those commodities typically transported via intermodal (*31*) and only those shipments with a distance of 500 miles or greater (*36*). The average load for a container/trailer is used for conversion, and the average train length in terms of TOFC/COFC count (*31*) is used to determine the number of intermodal trains equivalent to trucks hauled. The conversion methodology is as follows: (i) sort commodities transported by intermodal trains, (ii) convert tonnage of those commodities into equivalent container/trailer using average loading capacity, (iii) sum all container/trailer counts, and (iv) convert container/trailer counts to equivalent trainloads using average train length information. In intermodal transportation, truck haulage takes place from origin to delivery terminal and then from receiving terminal to destination. Therefore, every intermodal truck trip generates an empty truck trip. Thus, the number of container/trailer is doubled to obtain the intermodal truck flow. This procedure is carried out for all the demands that are transported via intermodal.

The conversion procedures were coded in Excel VBA to create freight O-D trip tables for truck, rail, and intermodal in 120 x 120 x 3 matrix form. It is assumed that road and rail infrastructure remains open for operation 365 days in a year. Using the aforementioned data sources and procedures, it is determined that in a single day in the base year (2007), there are 618,190 shipments transported by trucks, 1,415 shipments transported by trains, and 12,474 shipments transported via intermodal.

**Results and Discussions**

The solution algorithm was coded in MATLAB, and the experiments were run on a desktop computer with an Intel Core i7 3.40 GHz processor and 8 GB of RAM. The terms in the objective function are normalized to yield consistent units. This was accomplished by dividing the first term by the sum of truck demand and intermodal truck demand and second term by the sum of train demand and intermodal train demand. The stopping criterion used is the value of relative gap (change in value of objective function with respect to the value in previous iteration). The algorithm converged after 10 iterations in 686.50 seconds with a relative gap of $10^{-4}$. At convergence the value of the normalized objective function is 37.3594 hours. It should be noted that $\beta = 4$ is used here in the calculation of rail link delay.

The model was also solved using a classical algorithm (Frank-Wolfe). The Frank-Wolfe algorithm provides a normalized objective value of 37.3587 hours after 115 iterations and 2982.40 seconds of computational time. This result indicates that the Gradient-Projection algorithm is much more effective than the Frank-Wolfe algorithm in solving the proposed freight assignment model. This finding corroborates other studies which reported that the Gradient-Projection algorithm is superior to the Frank-Wolfe algorithm (e.g. *20*).

Among the four values tested for $\beta$, with $\beta = 2$ the flow on few links is very high, $\beta = 7$ the flow is reasonable, but the algorithm takes longer to converge, and $\beta = 15$ the flow results in very high travel time on some rail links. Therefore, for capturing freight train delay in the U.S. rail network, β = 4 is most suitable. Table 1 shows the percentage of link flow over capacity and link travel time of selected congested rail links, which were used to determine the best value for β. Note that, travel time is calculated based on the flow on corresponding link and flow on link opposite to it.

**TABLE 1 Comparison of β values**

| Link (Rail) | Percentage Increase in Flow over Capacity (Travel time in hour) | | | |
|---|---|---|---|---|
| Index | β = 2 | β = 4 | β = 7 | β = 15 |
| 76 | 43.3 (10.1) | 29 (7.9) | 27 (6.6) | 19.1 (39.9) |
| 81 | 36 (9.6) | 17 (9.5) | 11.7 (8.4) | 1.3 (11.5) |
| 268 | 27.3 (2.6) | 8.8 (2.2) | 2.1 (1.9) | 4.1 (2.5) |
| 279 | 29.5 (6.6) | 10.4 (6.1) | 5.6 (6.2) | 1 (7.2) |
| 392 | 30.3 (11.2) | 17.6 (8.8) | 4.6 (5.2) | 2.7 (15.5) |

The resulting user-equilibrium flow for the road network is shown in Figure 3(a) and for the rail network is shown in Figure 3(b). In Figure 3, the volume and spatial variation of freight traffic can be easily visualized by the thickness of the links.

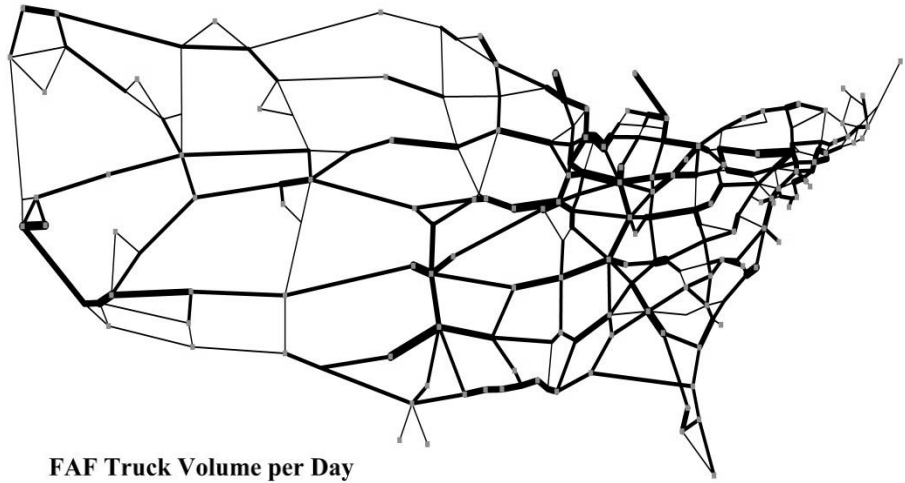

**FAF Truck Volume per Day**

— 0-5000   — 10000-15000   ▬ 20000-25000
— 5000-10000   ▬ 15000-20000

**(a)**

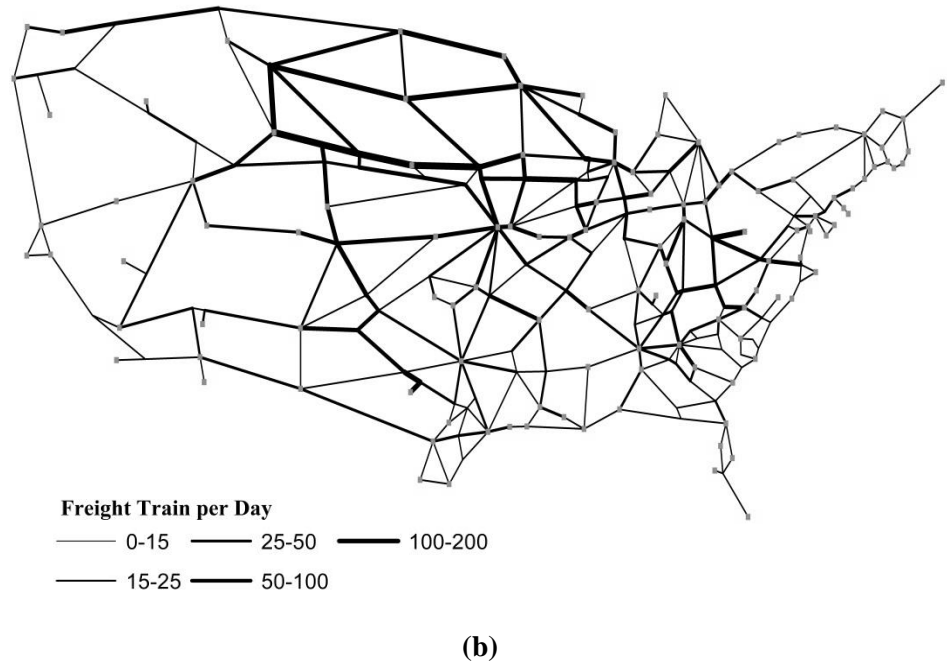

**Freight Train per Day**

—— 0-15   —— 25-50   —— 100-200

—— 15-25   —— 50-100

**(b)**

**FIGURE 3 Freight traffic assignment results: (a) truck in road network and (b) train in rail network.**

Figure 4(a) shows the FAF truck volume distribution on the U.S. national highway system for the year 2007. It shows truck flow patterns for trucks serving locations at least 50 miles apart and trucks not included in the "multiple modes and mail" (*37*). This truck flow pattern is very similar to the proposed model's projected user-equilibrium flow for the road network. Both maps indicate that there is high truck flow on interstates that traverse through California, Washington, Texas, Arkansas, Tennessee, Georgia, Florida, Michigan, Illinois, Indiana, New York, New Jersey, and Connecticut. This similarity suggests that the proposed model is capable of forecasting actual truck flows.

Figure 4(b) shows the 2005 freight trains and 2007 passenger trains volume per day in primary rail freight corridors in the U.S. (*31*). Though the proposed model's projected flow is only for freight train, this train flow pattern can be compared against the projected flow due to the fact that freight train volume far outnumbers passenger train volume in the U.S. The map indicates that there is high train flow on rail tracks that traverse through Washington, Montana, North Dakota, Arizona, New Mexico, Texas, Missouri, Wyoming, Nebraska, Iowa, Illinois, Indiana, Pennsylvania, Ohio, Georgia, New York, and New Jersey. The depicted train flow pattern and volume in most of the states are very similar to the proposed model's projected flow pattern. However, there exist a few discrepancies. The reason may be due to the difference in the demand between 2005 and 2007 and difference in methodology adopted to forecast freight flow. Note that Figure 4(b) is derived using annual survey data, whereas Figure 3(b) is derived from the equilibrium assignment procedure.

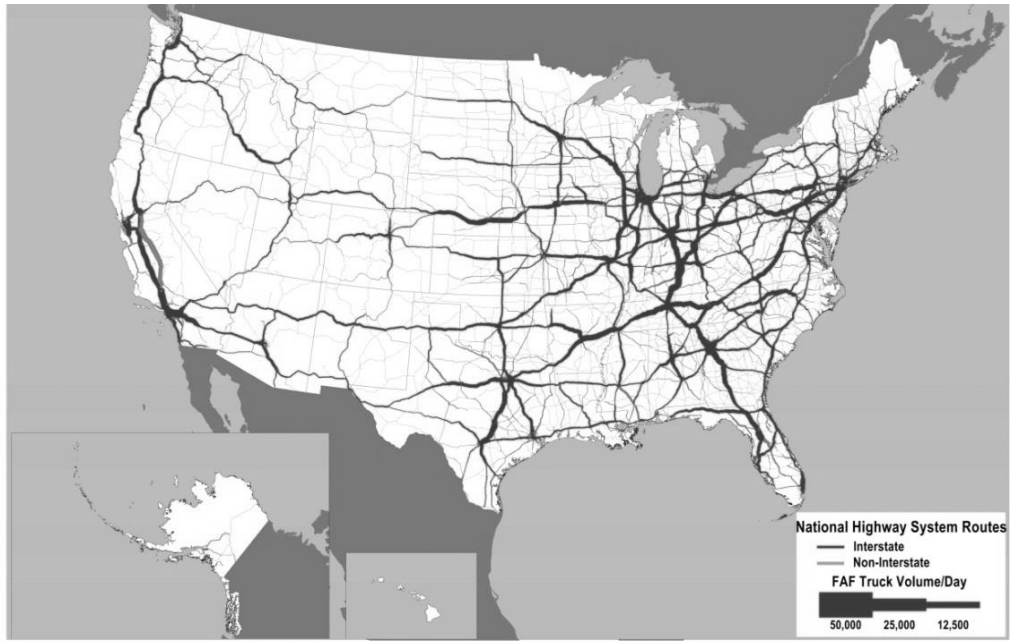

Note: Long-haul freight trucks typically serve locations at least 50 miles apart, excluding trucks that are used in movements by multiple modes and mail.
Source: U.S. Department of Transportation, Federal Highway Administration, Office of Freight Management and Operations, Freight Analysis Framework, version 3.4, 2012.

**(a)**

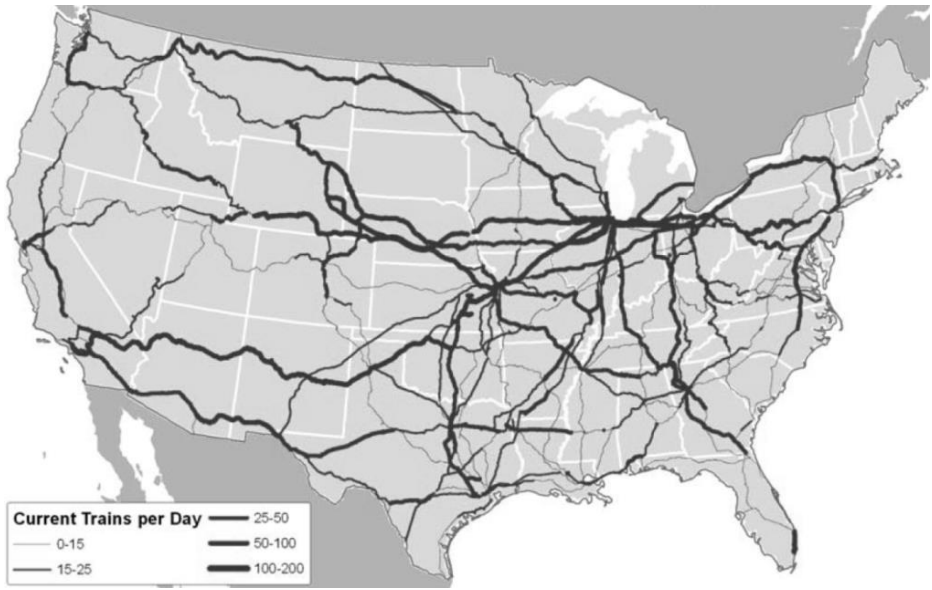

**(b)**

**FIGURE 4 Freight traffic volumes: (a) truck in U.S. highway system and (b) train in primary rail freight corridor.**

The proposed model's projected ton-miles are also compared quantitatively against those reported in the Commodity Flow Survey (CFS) and the FAF3. The results are reported in Table 2. In 2007, for the highway mode, the CFS reported freight ton-miles (*38*) is about 34% less than the FAF3 reported ton-miles. The difference in ton-miles between the proposed model and FAF3 and CFS is about 29% and 15%, respectively. Note that the FAF3 demand data was used as an input for the proposed model. Thus, the difference in ton-miles against FAF3 is reasonable because the proposed model only considered the contiguous U.S., and that it may have underestimated the intermodal demand. For the rail mode, the proposed model's projected ton-miles is about 14% less than that of the FAF3 data. This is reasonable for the same reasons mentioned previously. Overall, for both truck and rail demand, the proposed model appears to produce reasonable ton-miles value despite having a few simplifications, including a simplified network.

**TABLE 2 Freight Ton-Miles (Million) for Year 2007**

|              | FAF3      | CFS       | Proposed Model |
|--------------|-----------|-----------|----------------|
| Truck[a]     | 2,817,837 | 1,850,335 | 2,172,701      |
| Rail[b]      | 1,991,182 | 1,755,154 | 1,703,039      |

[a]Includes truck, and multiple modes and mail ; [b]Includes rail, and multiple modes and mail

**CONCLUSION**

This paper proposes a methodology for freight traffic assignment in large-scale road-rail intermodal networks. The proposed framework considers the lower level of a bi-level freight logistics problem, where the carriers' goals are to deliver the goods in a minimal amount of time. Given a set of freight demands between origins and destinations and designated modes (road-only, rail-only, and intermodal), the model finds the user-equilibrium freight flow. To obtain the solution for the model, a path-based algorithm based on the gradient projection algorithm is adopted. The proposed model is tested using the U.S. intermodal network and the FAF3 2007 freight shipment data. It is found that 4 is the most appropriate value for the β parameter when applying the Borndörfer et al. link performance function on the U.S. intermodal network. The results of the analysis, volume and spatial variation of freight traffic, show that the model produces equilibrium flow pattern that is very similar to the FAF3 flow assignment. The ton-miles values obtained from the model are also very close to those values reported in FAF3 and CFS. An attractive feature of the proposed model is that it converges within a few iterations and in about 11 minutes for a very large network. The model was also solved for the same network using the Frank-Wolfe algorithm, and results indicate that the Gradient-Projection algorithm is superior to the Frank-Wolfe algorithm in terms of convergence (i.e. fewer iterations) and computational time. The developed model can be used by transportation planners and decision makers to forecast freight flows and evaluate strategic network expansion options.

This study focuses on the freight assignment given the demands for various modes. In future work, we intend to build on this paper to formulate a bi-level mathematical program to capture the interactions between shippers' decisions and carriers' decisions. Also, we intend to expand the network to include additional modes (i.e., waterways, air).

**REFERENCES**

1. Hall, R. *Handbook of Transportation Science, 2ⁿᵈ edition*. Springer US, 2003.

2. Research and Innovative Technology Administration, Bureau of Transportation Statistics. *Freight Transportation: Global Highlights, 2010*. Publication BTS. U.S. Department of Transportation, Washington, D.C., 2010.

3. Strocko, E., M. Sprung, L. Nguyen, C. Rick, and J. Sedor. *Freight Facts and Figures 2013*. Publication FHWA-HOP-14-004. FHWA, U.S. Department of Transportation, 2013.

4. Brown, T. A., and A. B. Hatch. *The Value of Rail Intermodal to the U.S. Economy*. American Association of Railroads Report, 2002. www.intermodal.transportation.org/Documents/brown.pdf Accessed October 28, 2014.

5. Lowe, D. *Intermodal Freight Transport*. Elsevier Butterworth-Heinemann, 2005.

6. Chow, J. Y. J., S. G. Ritchie, and K. Jeong. Nonlinear Inverse Optimization for Parameter Estimation of Commodity-Vehicle-Decoupled Freight Assignment. *Transportation Research Part E: Logistics and Transportation Review,* Vol. 67, 2014, pp. 71-91.

7. Crainic, T.G., J. Ferland, and J. Rousseau. A Tactical Planning Model for Rail Freight Transportation. *Transportation Science,* Vol. 18, No. 2, 1984, pp. 165-184.

8. Guelat, J., M. Florian, and T. G. Crainic. A Multimode Multiproduct Network Assignment Model for Strategic Planning of Freight Flows. *Transportation Science,* Vol. 24, No. 1, 1990, pp 25-39.

9. Friesz, T. L., J. A. Gottfried, and E. K. Morlok. A Sequential Shipper–Carrier Network Model for Predicting Freight Flows. *Transportation Science,* Vol. 20, No. 2, 1986, pp. 80-91.

10. Fernández, J. E, J. De Cea, and R. Giesen. A Strategic Model of Freight Operations for Rail Transportation Systems. *Transportation Planning and Technology,* Vol. 27, No. 4, 2004, pp. 231-260.

11. Agrawal, B. B., and A. Ziliaskopoulos. Shipper–Carrier Dynamic Freight Assignment Model Using a Variational Inequality Approach. In *Transportation Research Record: Journal of the Transportation Research Board, No. 1966,* Transportation Research Board of the National Academies, Washington, D.C., 2006, pp. 60-70.

12. Loureiro, C. F. G., and B. Ralston. Investment Selection Model for Multicommodity Multimodal Transportation Networks. In *Transportation Research Record: Journal of the Transportation Research Board, No. 1522,* Transportation Research Board of the National Academies, Washington, D.C., 1996, pp. 38-46.

13. Kornhauser, A. L., and M. Bodden. Network Analysis of Highway and Intermodal Rail-Highway Freight Traffic. In *Transportation Research Record: Journal of the Transportation Research Board, No. 920,* Transportation Research Board of the National Academies, Washington, D.C., 1983, pp. 61-68.

14. Arnold, P., D. Peeters, and I. Thomas. Modelling a Rail/Road Intermodal Transportation System. *Transportation Research Part E: Logistics and Transportation Review,* Vol. 40, No. 3, 2004, pp. 255-270.

15. Mahmassani, H. S., K. Zhang, J. Dong, C. Lu, V. C. Arcot, and E. D. Miller-Hooks. Dynamic Network Simulation-Assignment Platform for Multiproduct Intermodal Freight Transportation Analysis. In *Transportation Research Record: Journal of the Transportation Research Board, No. 2032,* Transportation Research Board of the National Academies, Washington, D.C., 2007, pp 9-16.

16. Zhang, K., R. Nair, H. S. Mahmasanni, E. D. Miller-Hooks, V. C. Arcot, A. Kuo, J. Dong, and C. Lu. Application and Validation of Dynamic Freight Simulation- Assignment Model to Large-scale Intermodal Rail Network: Pan-European Case. In *Transportation Research Record: Journal of the Transportation Research Board, No. 2066,* Transportation Research Board of the National Academies, Washington, D.C., 2008, pp 9-20.

17. Yamada, T., B. F. Russ, J. Castro, and E. Taniguchi. Designing Multimodal Freight Transport Networks: A Heuristic Approach and Applications. *Transportation Science,* Vol. 43, No. 2, 2009, pp. 129-143.

18. Chang, T.-S. Best Routes Selection in International Intermodal Networks. *Computers & Operations Research,* Vol. 35, No. 9, 2008, pp. 2877-2891.

19. Hwang, T., and Y. Ouyang. Assignment of Freight Shipment Demand in Congested Rail Networks. In *Transportation Research Record: Journal of the Transportation Research Board,* Transportation Research Board of the National Academies, Washington, D.C., (forthcoming).

20. Jayakrishnan, R., W. Tsai, J. Prashker, and S. Rajadhyaksha. Faster Path-Based Algorithm for Traffic Assignment. In *Transportation Research Record: Journal of the Transportation Research Board, No. 1443,* Transportation Research Board of the National Academies, Washington, D.C., 1994, pp. 75-83.

21. Frank, M., and P. Wolfe. An Algorithm for Quadratic Programming. *Naval Research Logistics Quarterly,* Vol. 3, No. 1-2, 1956, pp. 95-110.

22. Chen, A., D. Lee, and R. Jayakrishnan. Computational Study of State-of-the-Art Path-Based Traffic Assignment Algorithms. *Mathematics and Computers in Simulation,* Vol. 59, No. 6, 2002, pp. 509-518.

23. Krueger, H. Parametric Modeling in Rail Capacity Planning. *Simulation Conference Proceedings*, Vol. 2, 1999, pp. 1194-1200.

24. Lai, Y., and C. Barkan. Enhanced Parametric Railway Capacity Evaluation Tool. In *Transportation Research Record: Journal of the Transportation Research Board, No. 2117,* Transportation Research Board of the National Academies, Washington, D.C., 2009, pp. 33-40.

25. Borndörfer, R., A. Fügenschuh, T. Klug, T. Schang, T. Schlechte, and H. Schülldorf. *The Freight Train Routing Problem*. ZIB Report, 2013. www.opus4.kobv.de/opus4-zib/frontdoor/index/index/docId/1899. Accessed May 5, 2014.

26. Southworth, F., D. Davidson, H. Hwang, B.E. Peterson, and S. Chin. *The Freight Analysis Framework Version 3 (FAF3): A Description of the FAF3 Regional Database and How It is Constructed*. FHWA, U.S. Department of Transportation, 2011.

27. Sheffi, Y. *Urban Transportation Networks: Equilibrium Analysis with Mathematical Programming Methods*. Prentice-Hall, Inc., New Jersey, 1985.

28. Bertsekas, D. P. On the Goldstein-Levitin-Polyak Gradient Projection Method. *Automatic Control, IEEE Transactions on,* Vol. 21, No. 2, 1976, pp. 174-184.

29. Center for Transportation Analysis (CTA), Oak Ridge National Laboratory (ORNL). *Intermodal Transportation Network*. www.cta.ornl.gov/transnet/Intermodal_Network.html. Accessed May 17, 2014.

30. Dowling,R., W. Kittelson, J. Zegeer, and A. Skabardonis. *Planning Techniques to Estimate Speeds and Service Volumes for Planning Applications*. Publication NCHRP Report 387. Transportation Research Board, Washington, D.C., 1997.

31. Cambridge Systematics, Inc. *National Rail Freight Infrastructure Capacity and Investment Study*. Association of American Railroads, 2007.

32. Standifer, G., and C. M. Walton. *Development of a GIS Model for Intermodal Freight*. U.S. Department of Transportation, 2000.

33. Freight Management and Operations (FMO), Federal Highway Administration (FHWA), U.S. Department of Transportation. *Freight Analysis Framework*. www.ops.fhwa.dot.gov/freight/freight_analysis/faf/index.htm. Accessed September 15, 2013.

34. Battelle. *FAF3 Freight Traffic Analysis*. Oak Ridge National Laboratory. www.faf.ornl.gov/fafweb/Data/Freight_Traffic_Analysis/faf_fta.pdf. Accessed December 22, 2013.

35. Hwang, T. Freight Transportation Demand Modeling and Logistics Planning for Evaluating Air Quality and Climate Impacts of Freight Systems. Ph.D. Dissertation, University of Illinois at Urbana-Champaign, Urbana, Illinois, 2013.

36. Slack, B. Intermodal Transportation in North America and the Development of Inland Load Centers. *The Professional Geographer,* Vol. 42, No. 1, 1990, pp. 72-83.

37. Freight Management and Operations (FMO), Federal Highway Administration (FHWA), U.S. Department of Transportation. *National Freight Transportation Maps*. www.ops.fhwa.dot.gov/freight/freight_analysis/nat_freight_stats/nhsavglhft2007.htm. Accessed May 5, 2014.

38. Margreta, M., C. Ford, and M. A. Dipo. *U.S. Freight on the Move: Highlights from the 2007 Commodity Flow Survey Preliminary Data*. Publication BTS-SR-018. RITA, U.S. Department of Transportation, 2009.